\documentclass[10pt]{amsart}
\usepackage{amssymb, amscd, psfrag, epsfig}

\newcommand{\Pro}{\mathbb{P}}
\newcommand{\C}{\mathbb{C}}
\newcommand{\Z}{\mathbb{Z}}
\newcommand{\rhobar}{{{\rho^*}}}

\newcommand{\usf}{\mathcal{U}^{SF}}
\newcommand{\vsf}{V^{SF}}

\newtheorem{theorem}{Theorem}[section]
\newtheorem{corollary}[theorem]{Corollary}
\newtheorem{definition}[theorem]{Definition}
\newtheorem{lemma}[theorem]{Lemma}
\newtheorem{proposition}[theorem]{Proposition}
\newtheorem{conjecture}[theorem]{Conjecture}
\newtheorem{question}[theorem]{Question}

\theoremstyle{remark} 
\newtheorem{remark}[theorem]{Remark}

\newcommand{\refT}[1]{Theorem~\ref{T:#1}}
\newcommand{\refC}[1]{Corollary~\ref{C:#1}}
\newcommand{\refP}[1]{Proposition~\ref{P:#1}}

\newcommand{\refL}[1]{Lemma~\ref{L:#1}}

\title{Bordism of semi-free $S^1$-actions}
\author{Dev Sinha}
\begin{document}

\begin{abstract}
We calculate geometric and homotopical  bordism rings 
associated to semi-free $S^1$ actions on complex manifolds, giving
explicit generators for the geometric theory.  The classification of 
semi-free actions with 
isolated fixed points up to cobordism complements similar
results from symplectic geometry.
\end{abstract}

\maketitle

\section{Introduction}

In this paper we describe both the geometric and homotopical bordism
rings associated to $S^1$-actions in which only the two simplest orbit 
types, namely fixed points and free orbits, are allowed.  Our work is of 
further interest in two different ways.  To make the computation 
of geometric semi-free bordism, in \refC{PTinj} 
we prove the semi-free case of what we call the 
geometric realization conjecture, 
which if true in general would determine the 
ring structure of geometric $S^1$-bordism from the ring structure of 
homotopical $S^1$-bordism given in \cite{Si00}.  
Additionally, we investigate semi-free actions with 
isolated fixed points as a first case, and that result 
is parallel to results from symplectic geometry.  
Let $\Pro(\C \oplus \rho)$ denote the space of 
complex lines in $\C \oplus \rho$ where $\rho$ is the standard 
complex representation of $S^1$ (in other words, the Riemann 
sphere with $S^1$ action given by the action of the unit complex numbers.)

\begin{theorem}\label{T:main} Let $S^1$ act semi-freely  with 
isolated fixed points on $M$,
compatible with a stable complex structure on $M$.  Then $M$ is
equivariantly cobordant to a disjoint union of products of
$\Pro(\C \oplus \rho)$.
\end{theorem}

This result should be compared with the second main result of 
\cite{TW00}, which states that when $M$ is connected
a semi-free Hamiltonian $S^1$ action 
on $M$ implies that $M$ has a perfect Morse function
which realizes the same Borel equivariant cohomology 
as a product of such $\Pro(\C \oplus \rho)$, as well as the same 
equivariant Chern classes. 
Our work also refines, in this case of isolated fixed points, results of
Stong \cite{Ston79}. 

As  \refT{main} led us to the more general 
computation of bordism of semi-free actions given in \refT{OmSF}, 
it would also be 
interesting to see if there is an analog of \refT{OmSF} for Hamiltonian
$S^1$-actions.  In general, the symplectic and 
cobordism approaches to transformation groups 
have considerable overlaps in language (for example, 
localization by inverting Euler classes of representations plays a key 
role in each theory), though the same words sometimes have different 
precise meanings.  A synthesis of these techniques might
be useful in addressing interesting questions within transformation groups
such as classifying semi-free actions with isolated fixed points.

In section~\ref{firstsec} of this paper we develop semi-free bordism theory
and give  a proof of \refT{main}.  We will see that the main ingredients are 
the Conner-Floyd-tom Dieck exact sequences, which are standard.  
In section~\ref{compsec} we compute semi-free
bordism theories.  In the final section, we review what is known
about $S^1$-bordism
and present a conjectural framework for the geometric theory.

The author would like to thank Jonathan Weitsman for 
stimulating conversations, and the referee, whose comments led to 
significant improvement of the paper.

\subsection{Notation}

If $X$ is a $G$-space, $X_+$ denotes $X$ with a disjoint basepoint 
with trivial action added.  If $V$ is a 
representation of $G$ equipped with a $G$-invariant inner product, 
let $S^V$ denote its  one-point compactification, let $D(V)$ be the unit 
disk in $V$, and let  $S(V)$ be
the boundary of $D(V)$, namely the unit sphere in $V$.
Let $\Omega^V X$ denote the space of based maps from $S^V$ 
to $X$, where $S^1$ acts by conjugation.
Let $X^{S^1}$ denote the fixed points of an $S^1$ action on $X$,
so that $Maps(X,Y)^{S^1}$ denotes the equivariant maps from
$X$ to $Y$.  Let $\bigoplus^n V = \bigoplus_{i=1}^n V$.
Let $\rho$ be 
the standard one-dimensional representation of $S^1$ and $\rhobar$  
its conjugate.

\section{First computations and \refT{main}}\label{firstsec}

The foundational results of this section are based on \cite{CCM96}, 
and the computational results parallel those of \cite{Si00}.  

\refT{main} follows from little more than the computation of Conner-Floyd
and tom Dieck exact sequences adapted for semi-free bordism.  
Because construction of these sequences 
is standard \cite{CF64, tD70, Alex72, Ston72, CCM96, Si01, Si01.2} , 
we will be brief in our exposition.

\begin{definition}
Let $\Omega^{SF}_*$ denote the bordism theory represented by 
stably complex (in the sense of Definition 28.3.1 of \cite{CCM96}) 
semi-free $S^1$-manifolds.  Bordisms between the manifolds must also
be semi-free (but see Remark~\ref{sfdef} below).  By 
equipping these manifolds and bordisms with equivariant maps to a space $X$ we
define  an equivariant homology theory $\Omega^{SF}_*(X)$. 
\end{definition}

Bordism theory is approachable in general because of its relation
to homotopy theory.  We choose a definition of homotopical equivariant
bordism with a relatively small amount of book keeping.

\begin{definition}
\begin{itemize}
\item Let $\vsf = \rho \oplus \C \oplus \rhobar$, with $S^1$-invariant
inner product defined 
through the standard inner products on $\rho$, $\C$ and $\rhobar$,
and let $\usf$ be 
$\bigoplus^\infty \vsf$.  Fix an isomorphism $\sigma : \usf \oplus \vsf \to
\usf$ sending $(w_1, w_2, \ldots)\oplus v$, with $w_i \in \vsf$, to $(v, w_1, w_2, \ldots)$.
\item Let $BU^{SF}(n)$ be the space of $n$-dimensional complex subspaces of 
$\usf$, topologized as the union over $k$
of $BU^{SF}(n,k)$, the $n$-dimensional 
subspaces of $\bigoplus^{n+k} (\rho \oplus \C \oplus \rhobar)$.
\item Let $\xi^{SF}(n)$ denote the total space of the
tautological bundle  over $BU^{SF}(n)$ with inner product inherited
from $\usf$, and let $TU^{SF}(n)$ be its
Thom space.  
\item Taking the direct sum of an $n$-dimensional subspace of $\usf$
with $\vsf$ defines a map 
$\xi^{SF}(n) \oplus \vsf$ to the total space of the tautological
bundle of $(n+3)$-dimensional subspaces of $\usf \oplus \vsf$, which 
through $\sigma$ is isomorphic to $\xi^{SF}(n+3)$.  Passing
to Thom spaces we get $\beta : S^{\vsf} \wedge TU^{SF}(n) \to TU^{SF}(n+3)$.
\item Let $MU^{SF}$ denote the $S^1$-spectrum with de-loopings
by semi-free representations built from the prespectrum $TU^{SF}$ whose
$V$th entry is $TU^{SF}(dim(V))$ and with bonding maps given by $\beta$.  
Explicitly, the $\bigoplus^k \vsf$th de-looping of the infinite loop space
associated to $MU^{SF}$ is given by the direct limit
$colim_n \; \Omega^{\bigoplus^{(n-k)} \vsf} TU^{SF}(3n)$, where the $\beta$
serve as maps in this directed system.
\end{itemize}
\end{definition}

Because any semi-free manifold can be embedded equivariantly in
some $\bigoplus^k \vsf$ (a direct application of transversality results
of \cite{Wass69} and the fact that $\rho$ and $\rhobar$ are the only 
representations which appear in the decomposition of the fiber of
the normal bundle to a fixed set), there is a Pontryagin-Thom
map from $\Omega^{SF}_*$ to $MU^{SF}_* = \pi_* MU^{SF}$.
We will see that this map is not an isomorphism but that nonetheless
$MU^{SF}_*$ is essential in studying $\Omega^{SF}_*$, 
in particular for proving \refT{main}.

\smallskip

The starting point in equivariant bordism is typically the use of
a filtration which can be traced back to Conner and Floyd \cite{CF64}.

\begin{definition}
\begin{itemize}
\item Define $i: MU_*(BS^1) \to \Omega^{SF}_*$ by taking a representative
$M$ mapping to $BS^1$ and pulling back 
the canonical $S^1$-bundle to get a principal $S^1$-bundle over $M$,
which is a free (and thus semi-free) $S^1$ manifold.  

\item  For a semi-free $S^1$-manifold $M$, the normal bundle of $M^{S^1}$ 
in $M$ will  have as the
representation type of the fiber a direct sum of $\rho$'s and $\rhobar$'s.  
Because $BU(n)$
classifies $n$-dimensional complex bundles, $MU_*(BU(n))$ is the bordism
module of stably complex manifolds with $n$-dimensional complex
bundles over them.  Let $F^{SF}_* = MU_*\left(\left(\bigsqcup_{n>0} 
BU(n)\right)^2\right)$ and define 
$\lambda : \Omega^{SF}_* \to F^{SF}_*$
as sending a semi-free bordism class to the bordism
class of the normal bundle of its fixed set, split according to appearance
of $\rho$ and $\rhobar$ in the fiber.
\item Define $\partial : F^{SF}_* {\to} MU_{*-2}(BS^1)$ as taking a 
manifold with a direct sum of
two bundles over it (classified by
maps to $BU(i) \times BU(j)$) imposing $S^1$ action
as $\rho$ on the summand of the first factor and $\rhobar$ on the second,
imposing an equivariant Hermitian inner product,
and then taking the unit sphere bundle of that $S^1$-bundle.
\end{itemize}
\end{definition}

\begin{theorem} \label{T:cfex}
The following  sequence is exact:
$$
0 \to \Omega^{SF}_* 
\overset{\lambda}{\to} F^{SF}_*
 \overset{\partial}{\to} MU_{*}(BS^1) \to 0.
$$
\end{theorem}

Note here that gradings are not preserved in the standard sense.  The middle
module must be graded so that $M$ mapping to $BU(i) \times BU(j)$ has degree
$dim(M) + 2(i+j)$.  The map $\partial$ lowers degree by two.

\begin{proof}[Outline of proof]
The maps $i$, $\lambda$ and $\partial$ coincide with the maps in the  
families exact sequence for the family $\{S^1, 1\}$
consisting of $S^1$ and the trivial group (see chapter
15 of \cite{CCM96}, or \cite{Si01}).  Exactness is straightforward and
pleasant to verify.  We claim that $i$ is the zero map.  It is well-known
that  $BS^1 = \C\Pro^\infty$ and 
$MU_*(\C\Pro^\infty)$ is generated by bordism representatives
$\C\Pro^n$ with their standard inclusions in $\C\Pro^\infty$ (see for example
Lemma 2.14 of part 2 of \cite{Adam74}).  The principal
$S^1$-bundle over $\C\Pro^n$ is equivariantly diffeomorphic to 
$S(\bigoplus^n \rho)$.  But this class is zero in $\Omega^{SF}_*$ since it bounds
$D(\bigoplus^n \rho)$.
\end{proof}

\begin{remark}\label{sfdef}
If we let $\Omega^{SF!}_*$ denote the image of semi-free bordism in 
unrestricted $S^1$ bordism thus allowing arbitrary bordism between semi-free
representatives, we see that $\Omega^{SF!}_*$ also fits in the exact 
sequence of \refT{cfex}, and thus is isomorphic to $\Omega^{SF}_*$ by the
five-lemma.
\end{remark}

The space $\bigsqcup_{n>0} BU(n)$ has a product  which corresponds 
to Whitney sum of bundles, through $BU(n)$'s role as the classifying space for
complex vector bundles.  Thus
$MU_*\left(\left(\bigsqcup_{n>0} BU(n)\right)^2\right)$ is a ring which we
identify as follows.

\begin{definition}
Let $X_{n,\rho} \in MU_{2n}\left( BU(1) \times BU(0) \right)$ be 
represented by
$\Pro^n$ mapping to $BU(1)$ by classifying the tautological line bundle.  
Let
$X_{n,\rhobar} \in MU_{2n}\left(BU(0) \times BU(1)\right)$ be defined 
similarly.
\end{definition}

\begin{proposition}\label{P:compFSF}
$ F^{SF}_*  \cong MU_*[X_{n,\rho}, X_{n, \rhobar}],$ where $n \geq 0$.
\end{proposition}

The proof is standard, as in Lemma 4.14 of part two of \cite{Adam74}, using the 
collapse of the Atiyah-Hirzebruch spectral sequence and the corresponding
computation for homology.  

\begin{corollary} \label{C:even}
$\Omega^{SF}_*$ is a free $MU_*$-module concentrated in even degrees.
\end{corollary}

\begin{proof}
Looking at the exact sequence of \refT{cfex} we see that the middle 
and right terms are free modules over $MU_*$.  The map $\partial$ is a split 
surjection, with one splitting given by sending the class represented by
$\C\Pro^n \hookrightarrow \C\Pro^\infty$ to the class represented by $D(\bigoplus^{n+1}
\rho)$, as in the outline of proof of \refT{cfex}.  As a submodule of $F^{SF}_*$,
$\Omega^{SF}_*$ is complementary to
the image of this splitting, and thus is free.
\end{proof}

We give one important example of computation of the map 
$\lambda$.

\begin{proposition} \label{P:complam}
$\lambda \left( \Pro(\C^n \oplus \rho) \right) = X_{n-1, \rho} + 
X_{0, \rhobar}^{n}$.
\end{proposition}


\begin{proof}
We use homogeneous coordinates on $\Pro(\C^n \oplus \rho)$.  
There are two possible components
of the fixed sets.   The points whose last coordinate is
zero constitute a fixed $\Pro^{n-1}$, whose normal bundle is
the tautological line bundle over  which  each fiber
is isomorphic to $\rho$ as a representation of $S^1$.  This
manifold with (normal) bundle defines exactly $X_{n-1, \rho}$.
There is also a fixed point in which all of the first $n$ coordinates are 
zero, and its normal bundle is $\bigoplus^n \rhobar$.  This fixed set
contributes a summand of $X_{0, \rhobar}^n$. 
\end{proof}

Next we introduce the analogue of \refT{cfex} for $MU^{SF}_*$, 
essentially the tom Dieck exact sequence.
We first need to develop Euler classes, which
play important roles in equivariant bordism.  
Consider $BU^{SF}(1)$, whose fixed set is three disjoint
copies of $BU(1)$.  The tautological bundle over $BU^{SF}(1)$
has fibers over these three fixed sets of $\rho, \C$ and $\rhobar$.

\begin{definition}
Let $\iota_\rho$ be the inclusion of a fiber isomorphic to $\rho$ 
over a fixed point  in the tautological
bundle over $BU^{SF}(1)$, noting that all such inclusions are homotopic.  Let
$T(\iota_\rho)$ denote the induced map on Thom spaces, and let 
$e_\rho \in MU^{SF}_{-2}$
be the composite $S^0 \to S^\rho \overset{T(\iota_\rho)}{\to} TU^{SF}(1)$.  
Let $e_\rhobar$ be defined similarly.
\end{definition}

The class $e_\rho$, when viewed as a class in 
$MU_{SF}^2(pt.)$ serves as the Euler class of $\rho$, viewed as a vector
bundle over a point.  

Next, we need to develop the analogue of $F^{SF}_*$.
Let $\Phi^{SF}_* = MU_*[ (BU \times \Z)^2]$, where multiplication on  
$(BU \times \Z)^2$ is the product of the standard Whitney sum multiplication
on each factor of $BU$ and addition on each factor of $\Z$.
By inclusion of $\bigsqcup_{n>0} BU(n)$ in $BU \times \Z$ (which is a group 
completion map, though we will not need that here), $F^{SF}_*$ 
maps to $\Phi^{SF}_*$.  The analogue of \refP{compFSF} is that
$$ \Phi^{SF}_* \cong MU_*[X_{0, \rho}^{\pm1}, X_{0, \rhobar}^{\pm 1},   
                                                X_{n,\rho}, X_{n, \rhobar} | n \geq 1] ,$$
where $X_{i, \rho}$ and $X_{i, \rhobar}$ are the images of the 
classes of the same name under the map from $F^{SF}_*$.
In particular, $X_{0, \rho}$ and $X_{0, \rhobar}$ are the unit 
classes in $(BU \times 1) \times (BU \times 0)$ and 
$(BU \times 0) \times (BU \times 1)$ respectively.  

\begin{theorem} \label{T:tdex}
There is a short exact sequence:
$$ 0 \to MU^{SF}_* \overset{\lambda}{\to} \Phi^{SF}_* \to MU_{*-2}(BS^1) \to 0.$$ 
The exact sequence of \refT{cfex} maps naturally to  this exact
sequence through Pontryagin-Thom maps.  The Pontryagin-Thom map
is  the identity 
on $MU_*(BS^1)$.  On the middle terms, $X_{i, \rho}$ and $X_{i, \rhobar}$
map to classes with the same names.  Moreover, $\lambda(e_{\rho}) =X_{0, \rho}^{-1}$ 
and $\lambda(e_{\rhobar}) = X_{0, \rhobar}^{-1}$. 
\end{theorem}

\begin{proof}[Outline of proof]
The proof of this theorem parallels the main results of \cite{tD70} and 
section four of \cite{Si00}.
The sequence in question is  the $MU^{SF}_*$ long exact sequence associated
to the cofiber sequence $ES^1_+ \to S^0 \to \widetilde{ES^1}$.
The middle term is of course $MU^{SF}_*$.  By either Adams' transfer
argument \cite{Ad2} 
or the fact that transversality holds in the presence
of free $G$-manifolds, $\widetilde{MU}^{SF}_*(ES^1_+)$ is isomorphic to
$MU_{*-1}(BS^1)$.  The map from  $MU_*(BS^1)$
to $MU^{SF}_*$ is zero since it factors through 
$i : MU_*(BS^1) \to \Omega^{SF}_*$,
which was shown to be zero in \refT{cfex}, so this long exact sequence
splits into short exact sequences. 

To identify $MU^{SF}_*( \widetilde{ES^1})$ as $\Phi^{SF}_*$ is a longer 
exercise.  The basic fact one uses is that if $X$ is semi-free and $Y$ is  contractible
when forgetting $S^1$-action (and both are CW-complexes)
then $Maps(X,Y)^{S^1}$ is homotopy equivalent
to $Maps(X^{S^1}, Y^{S^1})$ through the restriction map, since the fibers of 
this restriction map are spaces of (non-equivariant) maps into $Y$.
In analyzing $MU^{SF}_*( \widetilde{ES^1})$ one applies this fact to 
$Maps(S^V, \widetilde{ES^1}  \wedge TU^{SF}(n))$ to reduce to computing the
fixed sets  of these Thom spaces.   The fixed set $(TU^{SF}(n))^{S^1}$ 
is $\bigvee_{i+j+k = n} TU(i) \wedge (BU(j) \times BU(k))_+$ (see
Lemma~4.7 of \cite{Si00}).    Careful book keeping of the passage to spectra 
leads to the identification $$MU^{SF} \wedge \widetilde{ES^1}
\simeq \bigvee_{(i,j) \in \Z \times \Z} \Sigma^{2(i+j)} MU \wedge (BU \times BU)_+,$$ 
from which the isomorphism $MU_*(\widetilde{ES^1}) \cong \Phi^{SF}_*$ is immediate.


Identifying the Pontryagin-Thom map on the middle term
with the inclusion map from $F^{SF}_*$ to $\Phi^{SF}_*$ above is straightforward.
What remains is analysis of the Euler classs $e_\rho$ and $e_\rhobar.$ When
one passes to fixed sets, $e_\rho$ is represented by the inclusion 
$S^0 \hookrightarrow TU(0) \wedge (BU(1) \times BU(0))_+$. 
This class passes in the limit to the unit class in 
$MU_*((BU \times -1) \times (BU \times 0))$, which is the inverse of $X_{0, \rho}$.
The analysis of $e_{\rhobar}$ is similar. 

\end{proof}

In light of this theorem, we will usually express $\Phi^{SF}_*$ as
$MU_*[e_\rho^{\pm 1}, e_{\rhobar}^{\pm 1}, X_{n, \rho}, X_{n, \rhobar} | n \geq 1]$.
 From this theorem we deduce the following, whose first part is
an analogue of a theorem of Comeza\~na (28.5.4 of \cite{CCM96}) 
and L\"offler \cite{Lo72}.

\begin{corollary} \label{C:PTinj}
The Pontryagin-Thom map $\Omega^{SF}_* \to MU^{SF}_*$ is injective.
The following diagram from \refT{tdex} is a pullback square
$$
\begin{CD}
\Omega^{SF}_* @>\lambda>> F^{SF}_* \\
@V{P-T}VV @VVV \\
MU^{SF}_* @>>> \Phi^{SF}_*.
\end{CD}
$$
\end{corollary}

\begin{proof}
The horizontal maps are injective by Theorems~\ref{T:cfex} and~\ref{T:tdex},
the right vertical map is injective by inspection, so the left vertical map is 
injective by commutativity.

The horizontal maps have isomorphic cokernels, so the square is a pull-back
square through an elementary diagram chase.
\end{proof}

We use the phrase ``geometric realization'' to refer to the fact that this 
square is a pull-back, since it implies that any fixed-set data
which could be realized geometrically is so realized.
\refC{PTinj} will be the first ingredient 
in computing $\Omega^{SF}_*$ in the next section. 

Because homologically it is in negative degrees,
$e_\rho$ cannot be in the image of the Pontryagin-Thom map and thus might seem
exotic to the eyes of someone unfamiliar with equivariant bordism.  We will see
now that Euler classes can nonetheless be of great use in proving geometric 
theorems such as \refT{main}.

\begin{theorem}\label{T:main2} 
The intersection of $\lambda (MU^{SF}_*)$ with the
subring $\Z[e^{-1}_\rho, e^{-1}_\rhobar]$ is the subring $\Z[e^{-1}_\rho
+ e^{-1}_\rhobar]$.
\end{theorem}

Before proving this theorem, we deduce \refT{main} from it.

\begin{proof}[Proof of \refT{main}]
Let $M$ be a stably complex semi-free $S^1$-manifold with isolated
fixed points.  These isolated fixed points will have trivial normal bundles 
which are direct sums of $\rho$ and $\rhobar$.
Under $\lambda$, a fixed point with $\bigoplus^k \rho \oplus  \bigoplus^l \rhobar$ for
a normal bundle 
contributes $X_{0,\rho}^k X_{0, \rhobar}^l$.  By \refT{tdex},
this term maps to $e_\rho^{-k} e_\rhobar^{-l}$.
Therefore, $\lambda([M])$ lies in $\Z[e^{-1}_\rho, e^{-1}_\rhobar]$.

Applying \refT{main2}, $\lambda([M])$ lies in $\Z[e^{-1}_\rho
+ e^{-1}_\rhobar] $, which by \refP{complam} is
$\Z[\lambda\left(\Pro(\C \oplus \rho)\right)]$.
But by \refT{cfex} $\lambda$ is injective, so $[M]$ lies
in $\Z[\Pro(\C \oplus \rho)]$ in $MU^{SF}_*$.  Similarly,
by \refC{PTinj}, $[M]$ lies
in $\Z[\Pro(\C \oplus \rho)]$ in $\Omega^{SF}_*$, which
means that $M$ is equivariantly cobordant to a disjoint union
of products of $\Pro(\C \oplus \rho)$.
\end{proof}

Our main tool in the proof of \refT{main2} 
is  to use the augmentation map $\alpha : MU^{SF}_* \to MU_*$,
which takes a map $S^V \to TU^{SF}(n)$ and forgets the $S^1$ action.
Note that it is a map of rings.

\begin{proof}[Proof of \refT{main2}]
Let $R_*$ denote the subring $\Z[e^{-1}_\rho, e^{-1}_\rhobar]$ of $\Phi_*$.
Since $R_*$ is graded and lies in non-negative
degrees, we may proceed by induction on degree, 
focusing on homogeneous elements.  Suppose that 
$a_0 e_\rho^{-n} + a_1 e_\rho^{-(n-1)} e_\rhobar^{-1} +  \cdots + 
a_n e_\rhobar^{-n}$ is equal to $\lambda(x)$.  Consider 
$y = e_\rhobar (x - a_0 [\Pro(\C \oplus \rho)]^n)$.  The image
$\lambda(y)$ is in $R_*$ and is in degree $2(n-1)$, thus by inductive
hypothesis $y$ is in $\Z[\Pro(\C \oplus \rho)]$.
Hence $y = k [\Pro(\C \oplus \rho)]^{n-1}$ for some $k \in \Z$.  
Apply the augmentation map $\alpha$ to this equality. 
The image of $e_\rhobar$ under $\alpha$ is zero since $MU_{-2} = 0$, 
thus so is the image of $y$.
It is  well-known that $(\Pro^1)^{n-1}$ is non-zero in $MU_*$ for
any $n>0$, so $k$ must be zero.  This implies $y=0$, or since 
$e_\rhobar$ is not a zero divisor, $x =  a_0 [\Pro(\C \oplus \rho)]^n$.
The base case of this induction in degree zero is immediate since both 
$R$ and $\Z[e^{-1}_\rho + e^{-1}_\rhobar]$ consist only of the integers 
in that degree.
\end{proof}

\section{Computation of semi-free bordism}\label{compsec}

We turn our attention to homotopical semi-free bordism, following the
example of \cite{Si00}.  Let 
$Z_{n,\rho} \in \Omega^{SF}_*$ be
$[\Pro(\C^n \oplus \rho)]$, and by abuse let it also denote the image of this 
class under $\lambda$, which is equal to $X_{n-1, \rho} + X_{0, \rhobar}^{n}$
by \refP{complam}.  By further
abuse, let $Z_{n,\rho}$ also denote its image under the Pontryagin-Thom
map in $MU^{SF}_*$ as well as its image in $\Phi^{SF}_*$,
namely $X_{n-1, \rho} + e_\rho^{-n}$.   Let $Z_{n,\rhobar}$ be
defined (everywhere) similarly.  We may use $Z_{n, \rho}$ and $Z_{n, \rhobar}$
as generators of $F^{SF}_*$ and $\Phi^{SF}_*$.   By \refT{tdex} we have
the following.

\begin{proposition} \label{P:incl}
There is a sequence of inclusions
$$ MU_*[e_\rho, e_\rhobar, Z_{n,\rho}, Z_{n, \rhobar} | n \geq 2] \subset MU^{SF}_*
\subset MU_*[e_\rho^{\pm 1}, e_\rhobar^{\pm 1}, Z_{n,\rho}, Z_{n, \rhobar} | n \geq 2].$$
\end{proposition}

Thus, to understand $MU^{SF}_*$ is to understand divisibility by Euler classes,
which is traditionally done as part of a Gysin sequence.  
Recall $\alpha : MU^{SF}_* \to MU_*$, the augmentation map which forgets
$S^1$ action.

\begin{theorem} \label{T:eulaug}
The sequences  $0 \to MU^{SF}_{*+2} \overset{\cdot e_V}{\to} MU^{SF}_{*} \overset{\alpha}{\to}
MU_* \to 0,$
where $V$ is either $\rho$ or $\rhobar$, are exact.
\end{theorem}

\begin{proof}
Apply $\widetilde{MU}_{SF}^*$ to the cofiber sequence 
$S(\rho)_+ \overset{i}{\to} S^0 \overset{j}{\to} S^\rho$.  The middle term is
by definition
$MU_{SF}^*$.  Since $S(\rho)$ is a copy of the group 
$S^1$, an equivariant map is determined
by the image of one point so that
$Maps(S(\rho), X)^{S^1} = X$, for any $S^1$-space $X$ (with action
forgotten on the right-hand side), from which the similar statement follows
for spectra (see \cite{Ad2}) 
and in particular $MU^{SF}$.  The map $i^*$ is thus the augmentation map.

The identification of the remaining term is through a Thom isomorphism for 
$S^\rho$. Note that if an equivariant cohomology theory has such Thom
isomorphisms for all $S^V$ with $V$ complex it is said to be complex stable.
We roughly follow the construction of Thom isomorphisms for unrestricted
homotopical bordism given in section 10 of \cite{GrMa97}.
Unraveling definitions, we want to show that 
\begin{equation}\label{equality}
{\rm colim}_k \; \Omega^{\rho \oplus \bigoplus^k \vsf} TU^{SF}(3k) \simeq
{\rm colim}_k \; \Omega^{\C \oplus \bigoplus^k \vsf} TU^{SF}(3k).
\end{equation}

We start by choosing linear isomorphisms.  Choose coordinates
on $\usf = \bigoplus^\infty \vsf$ as 
$\bigoplus_{i=1}^\infty(v^\rho_i, v^\C_i, v^\rhobar_i)$ where $v^\rho_i \in \rho$,
$v^\C_i \in \C$ and $v^\rhobar_i \in \rhobar$.   Recall 
$\sigma : \usf \oplus \vsf \to \usf$, chosen to define
bonding maps for $MU^{SF}_*$, which in this notation sends
$\bigoplus_{i=1}^\infty(v^\rho_i, v^\C_i, v^\rhobar_i) \oplus (u^\rho, u^\C, u^\rhobar)$
to $\bigoplus_{i=1}^\infty(w^\rho_i, w^\C_i, w^\rhobar_i)$, where $w^\rho_1 = u^\rho$
and $w^\rho_i = v^\rho_{i-1}$ for $i>1$.  The vectors $w^\C_i$ and $w^\rhobar_i$
are defined similarly.  
Define 
$\sigma_1 : \usf \oplus \rhobar \oplus \bigoplus^2 \C \overset{\cong}{\longrightarrow} \usf$
by
$$\bigoplus_{i=1}^\infty(v^\rho_i, v^\C_i, v^\rhobar_i) \oplus u^\rhobar
\oplus u^\C_1 \oplus u^\C_2 \mapsto 
\bigoplus_{i=1}^\infty(w^\rho_i, w^\C_i, w^\rhobar_i),$$
where 

$$
w^\rho_i = v^\rho_i  \hspace{0.9in}
   w^\C_i = \begin{cases}  u^\C_i \;\; i \leq 2 \\  v^\C_{i-2} \;\; i> 2 \end{cases}  \hspace{0.9in}
        w^\rhobar_i = \begin{cases}  u^\rhobar_i \;\; i =1\\  v^\rhobar_{i-1} \;\; i> 1. \end{cases}
$$
Define $\sigma_2 : \usf \oplus \rhobar \oplus \bigoplus^2 \rho \overset{\cong}{\longrightarrow} \usf$ 
analogously so that the following diagram, in which the leftmost arrows are
the obvious isomorphisms which reorder coordinates, commutes:
\psfrag{A}{$\usf \oplus \bigoplus^2 \vsf$}
\psfrag{B}{$\usf \oplus (\rhobar \oplus \bigoplus^2 \C)$}
\psfrag{O}{$ \oplus (\rhobar \bigoplus^2 \rho) $}
\psfrag{I}{$\sigma_1 \oplus id$}
\psfrag{C}{$\usf \oplus  (\rhobar \oplus \bigoplus^2 \C)$}
\psfrag{L}{$\sigma_2$}
\psfrag{J}{$\sigma \oplus id$}
\psfrag{E}{$\usf \oplus \vsf$}
\psfrag{M}{$\sigma$}
\psfrag{H}{$\usf$.}
\psfrag{F}{$\usf \oplus (\rhobar \bigoplus^2 \rho) $}
\psfrag{Q}{$\oplus  (\rhobar \oplus \bigoplus^2 \C)$}
\psfrag{K}{$\sigma_2 \oplus id$}
\psfrag{G}{$\usf \oplus (\rhobar \oplus \bigoplus^2 \rho) $}
\psfrag{N}{$\sigma_1$}

$$\includegraphics[width=5.8in]{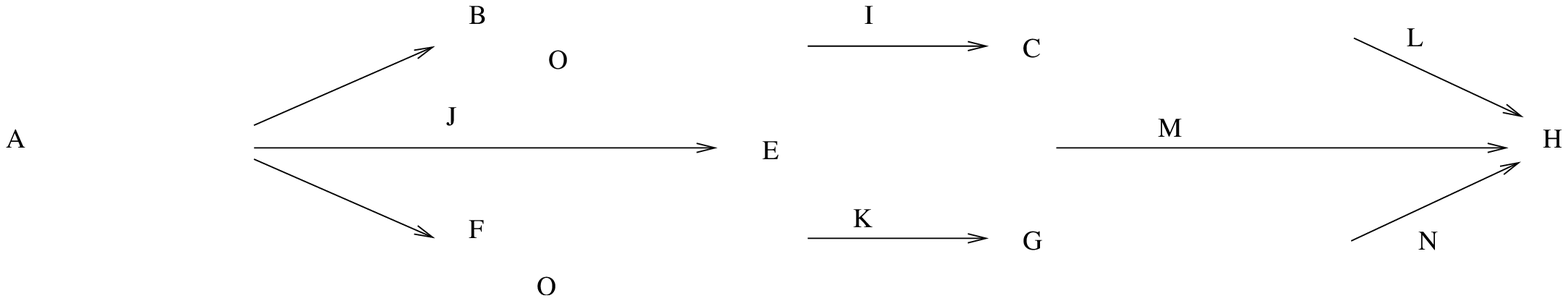}$$

On passage to Thom spaces, $\sigma_1$ defines a map 
$S^{(\rhobar \oplus \bigoplus^2 \C)} \wedge TU^{SF}(3k)$ to 
$TU^{SF}(3k+3)$, or by adjointness $T(\sigma_1)^\dagger : TU^{SF}_{3k} 
\to \Omega^{(\rhobar \oplus \bigoplus^2 \C)} TU^{SF}(3k+3)$.  Define
$$\beta_1 : \Omega^{(\rho \oplus \bigoplus^k \vsf)} TU^{SF}(3k)
\to \Omega^{(\C \oplus \bigoplus^{k+1} \vsf)} TU^{SF}(3k+3)$$
as sending an $f : S^{(\rho \oplus \bigoplus^k \vsf)} \to TU^{SF}(3k)$
to its composite 
with $T(\sigma_1)^\dagger$, using adjointness and the standard isomorphism
$(\rho \oplus \bigoplus^k \vsf) \oplus (\rhobar \oplus \bigoplus^2 \C) \cong
\C \oplus \bigoplus^{k+1} \vsf$ to get an element of the range.

Define $\beta_2 : \Omega^{(\C \oplus \bigoplus^k \vsf)} TU^{SF}(3k)
\to \Omega^{(\rho \oplus \bigoplus^{k+1} \vsf)} TU^{SF}(3k+3)$ similarly by using
$\sigma_2$.  By the commutativity of the diagram involving $\sigma$, $\sigma_1$
and $\sigma_2$ above, as well as standard facts about associativity of smash products
and adjointness, $\beta_2 \circ \beta_1$ and $\beta_1 \circ \beta_2$ 
coincide with the structure maps in the colimits of Equation~\ref{equality}, 
so that $\beta_1$ and $\beta_2$
give rise to the isomorphism of Equation~\ref{equality}.
Moreover, the map defining $j^*$ at the prespectrum level
composed with $\beta_1$ coincides with the definition of $e_\rho$,
so that $j^*$ is multiplication by $e_\rho$.

Finally, $MU_*$ is concentrated in even degrees, as is
$MU^{SF}_*$ since by \refT{tdex} it is a sub-algebra of $\Phi^{SF}_*$ which is so.
Therefore this long exact sequence
breaks up into short exact sequences, as stated.
\end{proof}

We introduce operations in $MU^{SF}_*$ which are essentially division
by Euler classes.  We will see below that these operations have a
geometric representation.

\begin{definition}
\begin{itemize}
\item Let $\sigma$ be the canonical (up to homotopy) splitting of the
augmentation map
$\alpha$, defined by taking some $S^m \to TU(n)$, suspending
it by $\bigoplus^m (\rho \oplus \rhobar)$ to get a map from $S^{\bigoplus^m \vsf}$ to
a Thom space which is chosen as a subspace of $TU^{SF}(n)$. 

\item Let $\Gamma_\rho : MU^{SF}_* \to MU^{SF}_{*+2}$  
(respectively $\Gamma_\rhobar$)
be the splitting of multiplication
by $e_\rho$ (respectively $e_\rhobar$)
which arises from the canonical splitting 
of $\alpha$ through \refT{eulaug}.

\item If $I$ is a sequence of $\rho$ and $\rhobar$, let $\Gamma_I(x)$ be
the composite of the corresponding $\Gamma_\rho$ and $\Gamma_\rhobar$ 
applied to $x$.  For example, $\Gamma_{\rho \rhobar}(x) = \Gamma_\rho
\Gamma_\rhobar (x)$.

\item For $x \in MU^{SF}_*$ let $\overline{x} = \sigma \circ \alpha(x)$.
\end{itemize}
\end{definition}

The following lemma is immediate from the fact that $e_\rho \Gamma_\rho(x)
= x - \overline{x}$.

\begin{lemma}\label{L:gamlam}
$\lambda(\Gamma_\rho(x))
= e_\rho^{-1} (\lambda(x) - \alpha(x))$ and similarly
$\lambda(\Gamma_\rhobar(x)) = e_\rhobar^{-1} (\lambda(x) - \alpha(x))$.
\end{lemma}

We are now ready for our first computation.

\begin{definition}
Let $B$ be the set of $MU^T_*$ elements $\{ e_\rho, e_\rhobar, 
Z_{n, \rho},$ and $Z_{n, \rhobar}\}$ where $n 
\geq 2$.  Order $B$ as follows 
$$e_\rho < e_\rhobar < Z_{2, \rho} < Z_{2,\rhobar} 
            < Z_{3,\rho} < Z_{3,\rhobar} < \cdots.$$
\end{definition}

\begin{theorem}\label{T:MUSF}
$MU^{SF}_*$ is generated as a ring by classes 
$\Gamma_\rho^i \Gamma_\rhobar^j (x)$ where $x \in B$. Relations are
\begin{enumerate}
\item $e_\rho \Gamma_\rho(x) = x - \bar{x} = e_\rhobar \Gamma_\rhobar(x),$ \label{defgam}
\item $\Gamma_V(x) (y - \bar{y}) = (x - \bar{x}) \Gamma_V(y)$, where $V$ 
is $\rho$ or $\rhobar$, \label{mult}
\item $\Gamma_V (e_V \cdot x) = x$, where $V$ is $\rho$ or $\rhobar$. \label{gamev}
\item $\Gamma_{\rhobar} \Gamma_\rho(x) = \Gamma_{\rho} \Gamma_\rhobar(x) + 
\overline{\Gamma_\rho(x)}  \Gamma_{\rho}\Gamma_\rhobar(e_\rho)$, \label{move}
\item $\overline{e_V} = 0$. 
\end{enumerate}

$MU^{SF}_*$ is free as a module over $MU_*$ with an 
additive basis given by monomials $\Gamma_\rho^i
\Gamma_\rhobar^j(x) m$, where $x \in B$, $m$ is a monomial in the 
$y \geq x$ in $B$ and with the following restrictions: if $x = e_\rhobar$ then $j=0$;
if $x = e_\rho$ and $j\neq 0$ then no positive power of $e_\rhobar$ occurs
in $m$;  if  $i \neq 0$ then $j \neq 0$ and no positive power of $e_\rho$ occurs in $m$.
\end{theorem}

\begin{proof}
\refP{incl} implies that if 
$y \in MU^{SF}_*$ then for some $i$ and $j$, the product $ x = e_\rho^i e_\rhobar^j y$
is in the subalgebra of $MU^{SF}_*$ generated by $B$.  
Then  $\Gamma_\rho^i \Gamma_\rhobar^j (x) = y$.
By linearity, $\Gamma_\rho^i
\Gamma_\rhobar^j (x)$ is a sum of $\Gamma _\rho^i \Gamma_\rhobar^j (m)$
for some monomials $m$ in $B$.   There is a product formula
$$\Gamma_\rho(wz) = \Gamma_\rho(w) z
+ \overline{w} \Gamma_\rho(z),$$ and similarly for $\Gamma_\rhobar$,
as can be verified by applying $\lambda$,
which is injective, to both sides using \refL{gamlam}.
Thus, $\Gamma _\rho^i \Gamma_\rhobar^j (m)$ is a sum of products of 
$\Gamma _\rho^i \Gamma_\rhobar^j (b)$ for $b \in B$, which means
these classes generate.

Except for relation~\ref{move}, verification of
the relations is straightforward.  In each case one checks the
equality after  $\lambda$, which is injective, using \refL{gamlam} as needed.
For example, for relation~\ref{mult}, the image of both sides under
$\lambda$ is $e_V^{-1}(x - \bar{x})(y - \bar{y})$.  For relation~\ref{move}
we also need that $\overline{\Gamma_\rho(x)} = - \overline{\Gamma_\rho(x)}$,
which we derive as follows.  Take relation~\ref{defgam} that 
$x = e_\rhobar \Gamma_\rhobar(x)$
and apply the product  formula with $w = e_\rhobar$ and $z = \Gamma_\rhobar(x)$
to get that $$\Gamma_\rho(x) = \Gamma_\rho(e_\rhobar) \Gamma_\rhobar(x),$$
noting that the second term in the product formula vanishes since $\overline{e_\rhobar} = 0$.
If we apply the augmentation map to both sides, $\overline{\Gamma_\rho(x)} = - \overline{\Gamma_\rho(x)}$ will follow from computing that 
$\overline{\Gamma_\rho(e_\rhobar)} = -1$.  
Represent $\Gamma_\rho(e_\rhobar)$ as the composite $S^{\rho} \to S^\rhobar
\to TU^{SF}(1)$, where the first map is through complex conjugation and the
second is the unit map, which includes $S^\rho$ as the Thom space of a fiber
of the tautological bundle.  This composite represents $-1$ when the 
$S^1$ action is forgotten.

To show that the members of the additive basis 
$\Gamma_\rho^i \Gamma_\rhobar^j (x) m$ are linearly independent
over $MU_*$ we apply $\lambda$, after which the verification is
straightforward by looking  at the leading terms $e_\rho^{-i} e_\rhobar^{-j} x m$.

To complete the proof we show that one can use the relations to reduce to the 
additive basis.
Consider a product $\Gamma_{I_1} (x_1) \Gamma_{I_2} (x_2) \ldots \Gamma_{I_k}
(x_k)$ where $x_1$ is minimal among the $x_i$ in order within $B$. 
We may use relation~\ref{mult}, rewritten as
$\Gamma_\rho(x) y =  x \Gamma_\rho(y) - \bar{x} \Gamma_\rho(y) + \bar{y}
\Gamma_\rho (x)$ (and similarly for $\rhobar$)
to perform a reduction.  Choose $y$ to be $\Gamma_{I_1}(x_1)$ and $x$ to be
$\Gamma_{I'_2}(x_2)$ where $I'_2$ is $I_2$ with the first $\rho$ or $\rhobar$
removed,  to decrease either the number of operations $\Gamma_V$
which are applied to non-minimal
generators, in the cases of  $x \Gamma_\rho(y)$ and  $\bar{x} \Gamma_\rho(y)$,
or the number of non-minimal generators, in the case of $\bar{y} \Gamma_\rho(x)$.  
Inductively, we reduce to a sum of $\Gamma_I(b) m$, where
$m$ is a monomial in $B$ and $b$ is less than any generator which appears
in $m$.  Finally, consider some
$\Gamma_{I_1 \rhobar \rho I_2} (b) m$.  We 
decrease the number of $\rho$ and $\rhobar$ which are out of order 
by applying relations~\ref{move} to get 
$\Gamma_{I_1 \rho  \rhobar I_2} (x) m +
\overline{\Gamma_{\rho I_2}(x)} \Gamma_{I_1 \rho \rhobar}(e_\rho) m$.
Note each of these monomials still has $\Gamma_V$ applied only to 
a minimal element of $B$.  Inductively, we reduce to monomials in which
$\Gamma_\rho$ is applied after $\Gamma_\rhobar$.

\end{proof}

We now turn our attention to $\Omega^{SF}_*$, adding to the
short list of geometric bordism theories which have been 
computed \cite{Alex72, Kosn76, Si01.2}.  By the geometric realization
\refC{PTinj} we can deduce the structure
of $\Omega^{SF}_*$ algebraically from \refT{MUSF} and understanding
of the localization map $\lambda$.  We choose, in addition, to find 
explicit geometric representatives.

We start by 
making geometric constructions of $\Gamma_\rho$ and 
$\Gamma_{\rhobar}$ on classes represented by 
manifolds.  These constructions follow ones made by Conner and Floyd.

\begin{definition}\label{D:gam}
Define $\gamma(M)$ for any stably complex $S^1$-manifold to be the
stably complex $S^1$-manifold
$$\gamma(M) = M \times_{S^1} S^3 \sqcup (-\overline{M}) \times 
\Pro (\C \oplus \rho),$$ 
where $S^3$ has the standard Hopf $S^1$-action
and the $S^1$-action on $M \times_{S^1} S^3$ is given by
\begin{equation}\label{E:dgam}
 \zeta \cdot \left[ m, z_1, z_2 \right] = \left[ \zeta \cdot
      m, z_1, \zeta z_2 \right].
\end{equation} 
Define ${\gamma}^*(M)$ similarly with the quotient of $M \times 
S^3$ by the $S^1$ action in which $\tau$ sends $m, (z_1, 
z_2)$ to $\tau m, (\tau z_1, \tau^{-1} z_2)$ and with induced $S^1$ 
action on the quotient given by 
\begin{equation}\label{E:dgamb}
 \zeta \cdot \left[ m, z_1, z_2 \right] = \left[ \zeta \cdot
      m, z_1, \zeta^{-1} z_2 \right].
\end{equation} 

\end{definition}

\begin{proposition}\label{P:model}
Let $M$ be  a stably complex $S^1$-manifold.  Then 
$\Gamma_\rho [M] = [\gamma(M)]$ and $\Gamma_\rhobar [M] =
[{\gamma}^*(M)]$ in $MU^{S^1}_*$.
\end{proposition}

\begin{proof}
By \refL{gamlam} and the injectivity of $\lambda$, it suffices to check 
the fixed sets and normal data of $\gamma(M)$ and $\bar{\gamma}(M)$.
One type of fixed points of $\gamma(M)$ are points 
$\left[ m, z_1, z_2 \right] $
such that $m$ is fixed in $M$ and $z_2 = 0$.  This fixed set is  diffeomorphic
to $M^G$, and its normal bundle is the normal bundle of $M^G$ in 
$M$ crossed with  the representation $\rho$.
Crossing with $\rho$ coincides with 
multiplying by $e_\rhobar$ in $F^{SF}_*$.
The second set of fixed points are $\left[ m, z_1, z_2 \right] $
such that $z_1 = 0$.  This set of fixed points is diffeomorphic to 
$M$, and its normal bundle is the trivial bundle $\rhobar$.

Hence, if $x  = \lambda([M])$, then the image of
$[\gamma(M)]$ is $x e_\rho^{-1} + \overline{M} e_{\rhobar}^{-1}$.
By subtracting the image of $\overline{M} \times \Pro (\C \oplus \rho)$
we obtain $x e_\rho^{-1} - \overline{M} e_{\rho}^{-1}$.
By \refL{gamlam}, this is $\lambda(\Gamma_\rho([M]))$.
The analysis is similar for ${\gamma}^*(M)$.
\end{proof}

The classes $\Gamma_I(Z_{n, \rho})$ and 
$\Gamma_I(Z_{n, \rhobar})$ can thus be realized geometrically, 
as $Z_{n,\rho}$ and $Z_{n, \rhobar}$  are represented
by linear actions on projective spaces.
Additionally we have the following.

\begin{lemma} \label{L:p1}
$\Gamma_{\rho \rhobar}(e_\rho) = \Pro(\C \oplus \rho)$.
\end{lemma}

\begin{proof}
The equality of these classes also follows from computation of their image 
under $\lambda$.   \refP{complam} states that
$\lambda\left(\Pro(\C \oplus \rho)\right) = e_\rho^{-1} + e_\rhobar^{-1}$. 
To show that this is also 
$\lambda(\Gamma_\rho \Gamma_\rhobar (e_\rho))$, by applying 
\refL{gamlam} twice it suffices to know that 
$\overline{\Gamma_\rho (e_\rhobar)} = -1$, which was shown in the proof of 
\refT{MUSF}.  
\end{proof}

Given the general complexities of equivariant bordism,
in particular for the geometric theories, $\Omega^{SF}_*$ has a 
remarkably simple form.

\begin{theorem}\label{T:OmSF}
$\Omega^{SF}_*$ is generated as an algebra
over $MU_*$ by classes $\gamma^i (\gamma^*)^j \Pro(\C^n \oplus \rho)$ for $n \geq 1$
and $\gamma^i (\gamma^*)^j \Pro(\C^n \oplus \rhobar)$ where $n \geq 2$.
Relations are
\begin{enumerate}
\item $\gamma(x) (y - \bar{y}) = (x - \bar{x}) \gamma(y)$, 
and similarly for $\gamma_*$, \label{multg}
\item $\gamma^* \gamma(x) = \gamma \gamma^*(x) + 
\overline{\gamma(x)}   \Pro(\C \oplus \rho)$, \label{moveg}
\end{enumerate}
where $x$ and $y$ can be any stably complex $S^1$-manifolds, in particular
those in the generating set above.
An additive basis is given by monomials $\gamma^i (\gamma^*)^j (x) m$ where $m$ 
is a monomial in $\Pro(\C^n \oplus \rho)$ and $\Pro(\C^n \oplus \rhobar)$ and $x$ is
such a projective space of smaller dimension than those appearing in $m$.

\end{theorem}

\begin{proof}
We start with \refC{PTinj}, which at the level of coefficients looks like
$$
\begin{CD}
\Omega^{SF}_* @>>> F^{SF}_* = 
     MU_*[e_\rho^{-1}, e_\rhobar^{-1}, Z_{n, \rho}, Z_{n,\rhobar} | n \geq 2] \\ 
@VVV @VVV \\
MU^{SF}_* @>>> \Phi^{SF}_* = 
    MU_*[e_\rho^{\pm 1}, e_\rhobar^{\pm 1}, Z_{n, \rho}, Z_{n,\rhobar} | n \geq 2],
\end{CD}
$$
where $n>0$.
All maps are inclusions, so we are looking to characterize the elements in
$MU^{SF}_*$ which map to $F^{SF}_*$.  Observe that $F^{SF}_*$
is an $MU_*$-direct summand  of $\Phi^{SF}_*$ .  A complementary submodule $C_*$ is
the submodule generated by reduced monomials 
in which a strictly positive power of $e_\rho$ or
$e_\rhobar$ appears.   We analyze the image under $\lambda$  of each additive
basis element from \refT{MUSF} in terms of the $F^{SF}_* \oplus C_*$ decomposition
of $\Phi^{SF}_*$.

Consider the basis element $y = \Gamma_\rho^i \Gamma_\rhobar^j(x) m$ in
which $x$ is an element of the generating set $B$ of \refT{MUSF} and
$m$ is a monomial in the elements of $B$, each greater than or equal to
$x$ in the ordering on $B$ and with additional provisions of $x = e_\rho$ or  $e_\rhobar$.  
This $y$ maps to $F^{SF}_*$ if $x = Z_{i, \rho}$ or
$Z_{i, \rhobar}$ because  by \refL{gamlam} $\lambda(\Gamma_I(Z_{i,\rho}))$ and 
$\lambda(\Gamma_I(Z_{i,\rhobar}))$ are polynomials
over $MU_*$ in $e_\rho^{-1}$, $e_\rhobar^{-1}$ and $Z_{i, \rho}$ or
respectively $Z_{i, \rhobar}$ and $m$ is a monomial
in $Z_{n, \rho}$ and $Z_{n, \rhobar}$ for some 
$n \geq i$ by the ordering on $B$.     Next we focus on when $x = e_\rho$.
By applying \refL{gamlam} we see that  
$\lambda(\Gamma_\rhobar^j(e_\rho) )= e_\rhobar^{-j} e_\rho + P$, where
$P$ is a polynomial in $e_\rhobar^{-1}$ over $MU_*$.   Continuing we see 
$\lambda(\Gamma_\rhobar^i \Gamma_\rhobar^j(e_\rho)) = e_\rho^{1-i} e_\rhobar^{-j} + Q$,
where $Q \in MU_*[e_\rho^{-1}, e_\rhobar^{-1}]$.
Recall that for the basis element
$y = \Gamma^i_\rho \Gamma^j_\rhobar (e_\rho) m$ with $i,j>0$, the generators
$e_\rho$ and  $e_\rhobar$ do not appear in $m$.   We deduce  that $\lambda(y)$
is in $F^{SF}_*$ since both 
$\lambda(\Gamma_\rhobar^i \Gamma_\rhobar^j(e_\rho))$ and $\lambda(m)$ are.

There are three classes of basis elements remaining, namely $e_\rho^i e_\rhobar^j m$ with
$i$ or $j > 0$, 
$\Gamma_\rhobar^i (e_\rho) e_\rho^j m$ with  $i>0$ and 
$\Gamma_\rho^i(e_\rhobar) e_\rhobar^j m$ with $i>0$,
where $m$ is a monomial in $MU_*[Z_{n, \rho}, Z_{n, \rhobar} | n \geq  2]$.  We take
the image under $\lambda$ and project onto $C_*$ to get  $e_\rho^i e_\rhobar^j m$,
$e_\rho^{j+1} e_\rhobar^{-i} m$ and $e_\rho^{-i} e_\rhobar^{j+1}m$ respectively.
These three kinds of classes 
are linearly independent taken all together
in $C_*$ (in fact, they form a basis as $m$ varies
over all possible monomials).  

Summarizing, we have shown that the additive basis
elements for $MU^{SF}_*$ fall into two groups, one group which maps to $F^{SF}_*$
and one group whose projections onto $C_*$ is linearly independent.  
Therefore, the only elements of $MU^{SF}_*$ which can map to $F^{SF}_*$
are in the span of the first group.  By \refC{PTinj} the first group serves as 
an additive basis for $\Omega^{SF}_*$.  

We will verify the additive
basis stated in the theorem only after we use the current additive basis to 
check that $\Omega^{SF}_*$ is generated as an algebra
by classes  $\gamma^i (\gamma^*)^j \Pro(\C^n \oplus \rho)$ and
$\gamma^i (\gamma^*)^j \Pro(\C^n \oplus \rhobar)$.  
By \refP{model}, $\gamma^i (\gamma^*)^j \Pro(\C^n \oplus \rho)$ represents 
$\Gamma_\rho^i \Gamma_\rhobar^j (Z_{n, \rho})$.   These generate
the additive basis elements of the form 
$\Gamma_\rho^i \Gamma_\rhobar^j(x) m$ where
$x = Z_{i, \rho}$ or $Z_{i, \rhobar}$.
To see that $\Gamma_\rho^i \Gamma_\rhobar^j(e_\rho)$ where $i,j > 0$
is in this subalgebra, first note that it is true for $i, j = 1$ by \refL{p1}.
We apply relation~\ref{move} from \refT{MUSF} to reduce to
this case as follows
\begin{align*}
\Gamma_\rho^i \Gamma_\rhobar^j(e_\rho)
  &=  \Gamma_\rho^{i-1} \Gamma_\rhobar \Gamma_\rho \Gamma_\rhobar^{j-1}(e_\rho)
        - \overline{\Gamma_\rho \Gamma_\rhobar^j(e_\rho)} \Gamma_\rho \Gamma_\rhobar (e_\rho)\\
   &= \cdots
      = \Gamma_\rho^{i-1} \Gamma_\rhobar^{j-1} \left( \Gamma_\rho \Gamma_\rhobar (e_\rho) \right) + Q,
\end{align*} 
where $Q \in MU_*[\Gamma_\rho \Gamma_\rhobar (e_\rho)]$.   
We see that $Q$ is in our subalgebra by \refL{p1}, 
which along with \refP{model} implies that 
$\Gamma_\rho^{i-1} \Gamma_\rhobar^{j-1} \left( \Gamma_\rho \Gamma_\rhobar (e_\rho) \right) = \gamma^{i-1} (\gamma^*)^{j-1} \Pro(\C \oplus \rho)$.  We deduce
that $\Gamma_\rho^i \Gamma_\rhobar^j(e_\rho)$
is in our subalgebra, so that all additive basis elements for $\Omega^{SF}_*$
are generated by the classes as stated.

The reduction to the additive basis given in the statement of the
theorem, and thus the proof that relations are complete, is similar to that given in \refT{MUSF}.
Given a monomial in $\gamma^i (\gamma^*)^j \Pro(\C^n \oplus \rho)$ and
$\gamma^i (\gamma^*)^j \Pro(\C^n \oplus \rhobar)$ we use relation~\ref{multg}
to reduce to monomials in which the operations $\gamma$ and $\gamma^*$ are applied
to only the projective space of the smallest dimension, 
and then use relation~\ref{moveg}
to reorder the operations.
\end{proof}

\section{Further directions in geometric bordism}

We are led to
ask about geometric bordism for unrestricted $S^1$ actions or
for actions by other groups.  
Bordism which is equivariant with respect to $\Z/p$ behaves similarly
to semi-free bordism, as expected.  
The Conner-Floyd and tom Dieck exact sequences are well-known
in those cases (indeed, it is for $\Z/p$ that these sequences first appeared in
\cite{CF64} and \cite{ tD70}),
and the theories were computed in \cite{Kosn76, Kriz99, Si99}, though the description
is complicated by the classes which are not restrictions from $\Omega^{SF}_*$.  
As in \refC{PTinj},  these theories fit in a pullback square
$$
\begin{CD}
\Omega^{U, \Z/p}_* @>>> F^{\Z/p}_*  \\
@VVV  @VVV \\
MU^{\Z/p}_* @>>> \Phi^{\Z/p}_*,
\end{CD}
$$
which follows because the kernels and cokernels of the horizontal 
maps are the even and odd degrees, respectively, of $MU_*(B\Z/p)$.
From this one can recover the Kosniowski generators from those of \cite{Si99}.
Note that  Kriz in \cite{Kriz99}  gave the first computation of $MU^{\Z/p}_*$,
but the relationship with the Kosniowski generators of geometric
bordism is not clear in Kriz's approach.

Less is known about $\Omega^{U, S^1}_*$, but we give a conjectural framework
as follows.  In \cite{Si00}, $MU^{S^1}_*$ was computed, and it has the following 
prominent features, much as we have seen for semi-free bordism:
\begin{itemize}
\item Basic classes include Euler classes $e_V$ and linear actions on projective spaces
$Z_{n,V} = [\Pro(\C^n \oplus V)]$ for all irreducible representations $V$.
\item There is a sequence of inclusions
$MU_*[e_V, Z_{n,V} | n \geq 2] \subset MU^{S^1}_* \subset \Phi_* = 
   MU_*[e_V^{\pm 1}, Z_{n,V} | n \geq 2]$,
where $V$ ranges over all irreducibles.
\item There are operations $\Gamma_V$ such that 
$e_V \Gamma_V(x) = x - \beta_V(x)$,
where $\beta_V(x)$ is restriction to $MU^{K(V)}_*$ followed 
by a splitting map back to $MU^{S^1}_*$.    Here $K(V)$ is 
the kernel of $V: S^1 \to \C^\times$.  Note that $\beta_V$ is not canonical
if $V = \rho$ or $\rhobar$.
\item $MU^{S^1}_*$ is generated over the operations $\Gamma^V$ by $e_V$
and $Z_{n, V}$.
\end{itemize}
There are also the following facts about the geometric theory:
\begin{itemize}
\item (Comeza\~na and L\"offler) The Pontryagin-Thom map $\Omega^{S^1}_* \to
MU^{S^1}_*$ is injective.
\item Under the inclusion $MU^{S^1}_* \to \Phi_*$,  the geometric theory
$\Omega^{U,S^1}_*$ maps to $F_* = MU_*[e_V^{-1}, Z_{n,V}]$.
\end{itemize}

A first important step towards understanding $\Omega^{U,S^1}_*$ would
be to establish the analogue of \refC{PTinj}, for which there are isolated
computations, as well as \refC{PTinj}, as evidence.

\begin{conjecture}
The square
$$
\begin{CD}
\Omega^{U, S^1}_* @>>> F_* \\
@VVV @VVV \\
MU^{S^1}_* @>>> \Phi_*
\end{CD}
$$
is a pull-back.
\end{conjecture}

This conjecture is likely to be approachable through the families filtration,
perhaps with $S^1$ replaced by $\Z/(p^2)$ as a starting point.
There would be two more steps needed to parallel our computation of
$\Omega^{SF}_*$.  

\begin{question}
Is there a version of the construction $\gamma$ for representations other
than $\rho$ and $\rhobar$?  In other words, given some $M$ can one find
a manifold which represents $\Gamma_V(M)$?
\end{question}

There is some doubt as to whether such a construction should even exist, 
given that embedded in such a construction would be a construction of
splitting maps $MU^{\Z/n}_* \to MU^{S^1}_*$, which are non-canonical
and chosen with some effort in \cite{Si00}.  A concrete starting point 
would be to search for a manifold whose fixed sets are  $D(\rho^2)$ 
crossed with the fixed sets  of 
$\Pro(\C^n \oplus \rho^3)$ and $\Pro(\C^n \oplus \rho)$ with its orientation
reversed.

We should add that even $\Gamma_\rho$ deserves more attention.
For example, what are the relationships between the equivariant
characteristic numbers (in both cohomology and $K$-theory) 
of $M$ and $\gamma(M)$?   How might $\Gamma_\rho$ be used to construct 
familiar classes in $MU_*$?   For example, in Proposition~6.5 
of \cite{Si00} we show
that $(\Gamma_\rho)^k(e_{\rho^n})$ form the coefficients of the
$n$-series.  

Finally, to compute $\Omega^{S^1}_*$ it would be helpful to 
understand the analogue of
\refT{main}, which promises to be much more difficult in the general setting.    
\refL{p1} that $\Gamma_\rho \Gamma_\rhobar (e_\rho) = [\Pro(\C \oplus \rho)]$
is surprising at first, since Euler classes seem unrelated to geometric ones.
But in fact all manifolds with framed fixed sets, in particular those with
isolated fixed sets, must arise within the description of $MU^{S^1}_*$ of 
\cite{Si00} as $\Gamma_I(x)$ where $x$ is a polynomial in $e_V$.
These constructions seem to be the most difficult part of describing
geometric classes within the homotopical setting, so once we proved
\refT{main} we knew \refT{OmSF} would be possible.  To provide a framework
for such investigation, we make the following.

\begin{conjecture}
Stably complex $S^1$ actions with isolated fixed points up to bordism are
generated by linear actions on projective spaces $\Pro(V_1 \oplus V_2 \oplus
\cdots \oplus V_k)$, where the weights of the $V_i$ are relatively prime.
\end{conjecture}

See Theorem~1.6 of \cite{Si00} for an example.  Taken all together, 
these questions and conjectures point to the following.

\begin{conjecture}
$\Omega^{S^1}$ is generated over geometric versions of the operations
$\Gamma_V$ by linear actions on projective spaces.
\end{conjecture}

\end{document}